\documentclass[leqno,12pt]{amsart}

\usepackage{amsmath,amstext,amssymb,amsopn,amsthm,mathrsfs}
\usepackage{enumerate}
\usepackage{color}
\usepackage{bbm}
\usepackage{verbatim} 
\usepackage{enumitem}
\usepackage{lipsum}

\newcommand\blfootnote[1]{%
  \begingroup
  \renewcommand\thefootnote{}\footnote{#1}%
  \addtocounter{footnote}{-1}%
  \endgroup
}

\numberwithin{equation}{section}

\allowdisplaybreaks

\DeclareMathOperator{\sgn}{sgn}
\DeclareMathOperator{\domain}{Dom}

\newtheorem{thm}{Theorem}[section]
\newtheorem{bigthm}{Theorem}

\newtheorem{lemma}[thm]{Lemma}
\newtheorem{corollary}[thm]{Corollary}

\newtheorem*{rem*}{Remark}

\newcommand{\DD}{\Delta\!\!\!\!\Delta}
\newcommand*\dd{{\rm d} }

\begin{document}

\blfootnote{


This research was supported by the National Science Centre of Poland within the project OPUS 2025/59/B/ST1/01786.
No AI/LLM was used to create this work.
}

\title[Maximal operator of the hybrid semigroup]
	{On the maximal operator of the hybrid
	\\ Hermite--Laguerre--Jacobi semigroup}

\author[A.\ Nowak]{Adam Nowak}
\address{Adam Nowak \newline
			Institute of Mathematics, Polish Academy of Sciences \newline
      \'Sniadeckich 8,
      00--656 Warszawa, Poland
      }
\email{adam.nowak@impan.pl}

\begin{abstract}
We consider the heat-diffusion semigroup related to expansions with respect to an orthogonal system
emerging by taking arbitrary finite tensor products of classical Hermite, Laguerre and Jacobi polynomials.
Our main result is the weak type $(1,1)$ estimate for the associated maximal operator.
This implies pointwise almost everywhere boundary convergence of the semigroup to initial values for all
initial conditions that are integrable against the orthogonality measure.
\end{abstract}

\maketitle

\section{Introduction} \label{sec:intro}

Harmonic analysis of orthogonal expansions is nowadays regarded as a classical part of harmonic analysis. It was initiated
by the seminal work of Muckenhoupt and E.M.\ Stein \cite{MuS}. This vast research area has been explored all over the world
in the last 60 years. A distinguished place in this theory belongs to expansions with respect to so-called
\emph{classical orthogonal polynomials}, namely \emph{Hermite}, \emph{Laguerre} and \emph{Jacobi polynomials} (see e.g. \cite{Sz}).
The measure spaces naturally associated with these expansions in dimension one are, respectively,
$$
\Big( \mathbb{R}, e^{-x^2}\, {\dd}x\Big), \qquad \Big( [0,\infty), x^{\alpha}e^{-x}\, {\dd}x \Big), \qquad
	\Big( [-1,1], (1-x)^{\alpha}(1+x)^{\beta} {\dd}x \Big);
$$
here $\alpha,\beta > -1$ are so-called parameters of type. Note that the first two are non-doubling measure spaces
while the last one is doubling (which means that concentric balls of comparable radii have comparable measures).
Multi-dimensional situations arise in a standard way, by taking products of one-dimensional
polynomials and measures of orthogonality.

In each of the three frameworks there is a naturally associated symmetric diffusion semigroup of operators $\{T_t\}_{t\ge 0}$:
Hermite (Ornstein--Uhlenbeck), Laguerre and Jacobi semigroup, respectively.
It is well known, see e.g.\ \cite{BGL}, that in dimension one and up to affine transformations and time rescaling,
these are the only conservative symmetric diffusion semigroups on an interval whose generators are second-order
differential operators, elliptic in the interior,
with a complete orthogonal basis of polynomial eigenfunctions of successive degrees.
Those semigroups are also transition semigroups for the Ornstein--Uhlenbeck, Laguerre and Jacobi stochastic diffusion processes,
which are of importance in theory and applications, in mathematics, physics and beyond.

Maximal operators
$$
T_*f = \sup_{t>0}|T_tf|
$$
of the three semigroups and their mapping properties are relevant since, for instance, they lead directly to boundary
convergence results for these semigroups as $t \to 0^+$.
Boundedness on $L^p$ spaces, $p>1$, of the three maximal operators follows from
Stein's maximal theorem \cite[p.\,73]{St}.
On the other hand, there is no $L^1$-boundedness. The natural substitute that still implies the boundary convergence
for $L^1$ functions is weak type $(1,1)$ boundedness, a property that is, however, much harder to prove.

For the Ornstein--Uhlenbeck semigroup maximal operator weak type $(1,1)$ is due to Muckenhoupt \cite{Mu} in dimension one,
and due to Sj\"ogren \cite{Sj} in higher dimensions. Later at least three other proofs of this result were found \cite{MRS,GMMST,AFS}.
In case of the Laguerre semigroup maximal operator the one-dimensional result was proved again by Muckenhoupt \cite{Mu},
while the multi-dimensional counterpart was claimed first in \cite{Di} and later in \cite{Sa2}, but we believe that the first
genuinely complete proof was delivered only in \cite{NSS}, by Sj\"ogren, Szarek and the author;
see the related comments in \cite[p.\,313]{NSS}.
For the Jacobi semigroup maximal operator weak type $(1,1)$ was shown by Sj\"ogren and the author \cite{NoSj}.
In that case the main obstacle was to find suitable bounds for the semigroup integral kernel, which were obtained
in \cite{CKP}, and independently in \cite{NoSj} under a slight restriction on the type parameters.
The result is also a straightforward consequence of the Gaussian upper bound for the multi-dimensional Jacobi heat kernel
in trigonometric parametrization, see the comments in \cite[p.\,349]{NSS}.

Now, one can think about mixed Hermite--Laguerre--Jacobi expansions, where the orthogonal system is defined by taking
joint tensor products of finite number of one-dimensional
Hermite, Laguerre and Jacobi polynomials (allowing for any polynomial type to be omitted from the product).
Then, after suitable re-numeration of axes, the corresponding semigroup of operators will have an integral kernel
which is a product of (multi-dimensional) integral kernels of Hermite, Laguerre and Jacobi semigroups.
The resulting \emph{hybrid Hermite--Laguerre--Jacobi semigroup} is again a symmetric diffusion semigroup, which includes as special
cases the three fundamental semigroups. In particular, boundedness on $L^p$, $p>1$, of the hybrid semigroup maximal operator
follows from Stein's maximal theorem \cite[p.\,73]{St}. But weak type $(1,1)$ does not follow from known results and is far from
obvious.

The aim of this paper is to prove weak type $(1,1)$ boundedness for the maximal operator of the
hybrid Hermite--Laguerre--Jacobi semigroup, see Theorem \ref{thm:main} in Section~\ref{sec:main}.
We believe this is an interesting result that, in particular, unifies the corresponding results in all
the three settings. Moreover, it contributes to the theory of maximal operators by
rather non-standard mixture of non-doubling and doubling situations,
a combination potentially creating technical challenges not present in each case separately.
Further originality comes from the fact that, to the best of our knowledge, hybrid Hermite--Laguerre--Jacobi
scenario has not yet appeared explicitly in the literature devoted to harmonic analysis of orthogonal expansions.

Our proof of Theorem~\ref{thm:main} relies on a gradual reduction of the problem to weak type $(1,1)$ estimate
for the multi-dimensional Laguerre semigroup maximal operator, which is known. This reveals that, in a sense, the Laguerre maximal
operator is superior to the Hermite and Jacobi ones, since it is more sophisticated and, moreover, it encompasses difficulties
occurring separately in analysis of the Hermite and Jacobi maximal operators.

The paper is organized as follows.
In Section~\ref{sec:prel} we briefly recall the Hermite, Laguerre and Jacobi settings. Then in Section~\ref{sec:mix} 
we introduce the mixed Hermite--Laguerre--Jacobi framework and the hybrid Hermite--Laguerre--Jacobi semigroup.
Section \ref{sec:main} contains statement of the main result (Theorem~\ref{thm:main}) and a corollary concerning
almost everywhere boundary convergence of the Hermite--Laguerre--Jacobi semigroup.
Finally, Section~\ref{sec:proof} is devoted to the proof of Theorem~\ref{thm:main}.

\subsection*{Notation}
Set $\mathbb{R}_{+} := [0,\infty)$.
Throughout the paper we denote by $H_n$, $L_n^{\alpha}$ and $P_n^{\alpha,\beta}$ the classical Hermite, Laguerre
and Jacobi polynomials, respectively, as defined e.g.\ in Szeg\H{o}'s monograph \cite{Sz}.
We use standard conventions concerning empty sums and products.
We write $X \lesssim Y$ to indicate that
$X \le C Y$ with a positive constant $C$ independent of significant quantities. We shall write
$X \simeq Y$ when simultaneously $X \lesssim Y$ and $Y \lesssim X$.

\section{Preliminaries} \label{sec:prel}

We now briefly describe the well-known one-dimensional Hermite, Laguerre and Jacobi settings.
Multi-dimensional frameworks arise in a standard way, by taking tensor products of one-dimensional situations, and they
are contained as special cases in the
general mixed Hermite--Laguerre--Jacobi setting that will be introduced subsequently.

\subsection{Hermite setting} \label{ssec:H}
Let
$$
\dd \gamma(x) := e^{-x^2}\, \dd x
$$
be the (non-normalized) Gaussian measure in $\mathbb{R}$. Hermite polynomials $H_{n}$, $n=0,1,2,\ldots$, form
an orthogonal basis in $L^2(\mathbb{R},\dd\gamma)$. They are eigenfunctions of the Ornstein--Uhlenbeck operator
$$
\DD_{H} := -\frac{1}2 \frac{\dd^2}{\dd x^2} + x \frac{\dd}{\dd x}, 
$$
one has
$$
	\DD_H H_n = n H_n, \qquad n \ge 0.
$$
The Ornstein--Uhlenbeck operator is symmetric and non-negative in $L^2(\mathbb{R},\dd \gamma)$
when considered initially on the subspace of polynomials.
Thus there is the natural self-adjoint extension in $L^2(\mathbb{R},\dd\gamma)$
of $\DD_H$ (denoted by the same symbol) acting initially on polynomials, which is non-negative and
whose spectral resolution is given in terms of Hermite polynomials.

The associated semigroup $\{\exp(-t\DD_H)\}_{t \ge 0}$ is called the the Ornstein--Uhlenbeck or the Hermite semigroup.
There is an integral representation
$$
\exp(-t \DD_H) f(x) = \int_{\mathbb{R}} K_t^{H}(x,y) f(y)\, \dd \gamma(y), \qquad x \in \mathbb{R}, \quad t >0,
$$
valid for\,\footnote{
\,The integral formula extends the action of the Hermite semigroup to $L^1(\dd\gamma)$ and thus to all
$L^p(\dd\gamma)$ spaces, $p \ge 1$, which are included in $L^1(\dd\gamma)$.
}
$f \in L^1(\dd\gamma)$, where the kernel is given by the explicit formula\,\footnote{
\,Computation of the series is due to F.G.\ Mehler, back in 1866; see e.g.\ \cite{Sj2} for the reference.
Therefore the kernel is called sometimes the Mehler kernel.
} 
\begin{align*}
K_t^H(x,y) & := \sum_{n \ge 0} e^{-tn} \frac{H_n(x)H_n(y)}{\|H_n\|^2_{L^2(\dd \gamma)}} \\
	& = \sqrt{ \frac{e^t}{2\pi \sinh t} }\, \exp\bigg( - \frac{e^{-t}(x^2+y^2)-2xy}{2\sinh t}\bigg),
		\qquad x,y \in \mathbb{R}, \quad t > 0.
\end{align*}
Note that $\|H_n\|^2_{L^2(\dd \gamma)} = \sqrt{\pi}\,2^n n!$.

\subsection{Laguerre setting} \label{ssec:L}
For a parameter $\alpha > -1$, define a gamma type measure in $\mathbb{R}_+$ by
$$
\dd \mu_{\alpha}(x) := x^{\alpha}e^{-x}\, \dd x.
$$
Laguerre polynomials $L_n^{\alpha}$, $n=0,1,2,\ldots$, constitute an orthogonal basis in $L^2(\mathbb{R}_+,\dd\mu_{\alpha})$
of eigenfunctions of the Laguerre operator\,\footnote{
\,In the definition of $\DD_{L,\alpha}$ the common factor $2$ is to make the Hermite and Laguerre semigroups `time' parameters $t$
consistent, which might be of importance for potential applications of the hybrid Hermite--Laguerre(--Jacobi) semigroup.
This adjustment is perhaps best explained in terms of the Ornstein--Uhlenbeck and Laguerre stochastic processes,
whose infinitesimal generators are a multi-dimensional variant of $\DD_H$ and $\DD_{L,\alpha}$, respectively.
Indeed, if $\{X_t : t \ge 0\}$ is a $d$-dimensional Ornstein--Uhlenbeck process, then the squared radial Ornstein--Uhlenbeck process
$\{|X_t|^2 : t \ge 0\}$ is a one-dimensional Laguerre process $\{Y_t^{\alpha} : t \ge 0\}$ with parameter $\alpha = d/2-1$.
See \cite[Appendix 1, Sections 24--26]{BS}. Without the factor $2$ the clocks of $|X_t|^2$ and $Y_t^{d/2-1}$ would run at
different speeds.
}
$$
\DD_{L,\alpha} := -2x\frac{\dd^2}{\dd x^2} - 2(\alpha+1-x) \frac{\dd}{\dd x};
$$
one has 
$$
	\DD_{L,\alpha} L_n^{\alpha} = 2 n L_n^{\alpha}, \qquad n \ge 0.
$$
The Laguerre operator is symmetric and non-negative in $L^2(\mathbb{R}_+,\dd \mu_{\alpha})$ when considered initially on the
subspace of polynomials.
The natural in this context self-adjoint extension in $L^2(\mathbb{R}_+,\dd\mu_{\alpha})$ of $\DD_{L,\alpha}$
is non-negative and its spectral resolution is given in terms of Laguerre polynomials.

The corresponding semigroup $\{\exp(-t\DD_{L,\alpha})\}_{t \ge 0}$ is called the Laguerre semigroup.
This semigroup admits an integral representation
$$
\exp(-t\DD_{L,\alpha})f(x) = \int_{\mathbb{R}_{+}} K_t^{L,\alpha}(x,y) f(y)\, \dd\mu_{\alpha}(y), \qquad x \in \mathbb{R}_+, \quad t>0,
$$
which is valid for $f \in L^1(\dd\mu_{\alpha})$. The integral kernel here is given explicitly by\,\footnote{
\,Computation of the series is possible thanks to the Hille-Hardy formula discovered in 1926, see e.g.\ \cite[p.\,220]{NoSj} for the
related comments and references.
} 
\begin{equation} \label{lag_form}
\begin{split}
& K_t^{L,\alpha}(x,y) := \sum_{n \ge 0} e^{-t2n} \frac{L_n^{\alpha}(x)L_n^{\alpha}(y)}{\|L_n^{\alpha}\|^2_{L^2(\dd\mu_{\alpha})}} \\
& = \frac{e^{t(\alpha+1)}}{2\sinh t} \exp\bigg( -\frac{e^{-t}}{2\sinh t}(x+y)\bigg) \big(\sqrt{xy}\,\big)^{-\alpha}
	I_{\alpha}\bigg( \frac{\sqrt{xy}}{\sinh t}\bigg), \qquad x,y \in \mathbb{R}_+, \quad t > 0;
\end{split}
\end{equation}
note that $\|L_n^{\alpha}\|^2_{L^2(\dd\mu_{\alpha})} = \Gamma(n+\alpha+1)/n!$.
In the above formula $I_{\alpha}$ is the modified Bessel function of the first kind and of order $\alpha$.
The values of the kernel when $x=0$ or $y =0$ exist, in the limiting sense, since
$I_{\alpha}(z) = z^{\alpha}/(2^{\alpha}\Gamma(\alpha+1)) + \mathcal{O}(z^{\alpha+2})$, $z \to 0^+$, cf.\ \cite[(5.7.1)]{Leb}.
For further reference, note that standard asymptotics for $I_{\alpha}$ with $\alpha > -1$
imply, cf.\ \cite[(5.16.5)]{Leb},
\begin{equation} \label{bes_est}
I_{\alpha}(z) \simeq z^{\alpha}(z+1)^{-\alpha-1/2}e^z, \qquad z > 0.
\end{equation}

\subsection{Jacobi setting}
Given parameters $\alpha, \beta > -1$, define a beta type measure in $[-1,1]$
$$
\dd\rho_{\alpha,\beta}(x) := (1-x)^{\alpha}(1+x)^{\beta}\, \dd x.
$$
Jacobi polynomials $P_n^{\alpha,\beta}$, $n=0,1,2,\ldots$, form an orthogonal basis in $L^2([-1,1],\dd\rho_{\alpha,\beta})$
of eigenfunctions of the Jacobi operator
$$
\DD_{J,\alpha,\beta} := -\big(1-x^2\big) \frac{\dd^2}{\dd x^2} - \big[\beta-\alpha-(\alpha+\beta+2)x\big] \frac{\dd}{\dd x};
$$
one has
$$
\DD_{J,\alpha,\beta} P_n^{\alpha,\beta} = n(n+\alpha+\beta+1) P_n^{\alpha,\beta}, \qquad n \ge 0.
$$
The Jacobi operator is symmetric and non-negative in $L^2([-1,1],\dd \rho_{\alpha,\beta})$ when considered initially on the
subspace of polynomials. The natural in this context self-adjoint extension in $L^2([-1,1],\dd\rho_{\alpha,\beta})$ of
$\DD_{J,\alpha,\beta}$ is non-negative and its spectral resolution is given in terms of Jacobi polynomials.

The semigroup $\{\exp(-t\DD_{J,\alpha,\beta})\}_{t\ge 0}$ is called the Jacobi semigroup. There is an integral representation
$$
\exp(-t\DD_{J,\alpha,\beta})f(x) = \int_{[-1,1]} K_t^{J,\alpha,\beta}(x,y)f(y)\, \dd\rho_{\alpha,\beta}(y),
	\qquad x \in [-1,1], \quad t > 0,
$$
valid for $f \in L^1(\dd\rho_{\alpha,\beta})$. Here the kernel is given by
$$
K_t^{J,\alpha,\beta}(x,y) := \sum_{n \ge 0} e^{-tn(n+\alpha+\beta+1)}
	\frac{P_n^{\alpha,\beta}(x)P_n^{\alpha,\beta}(y)}{\|P_n^{\alpha,\beta}\|^2_{L^2(\dd\rho_{\alpha,\beta})}};
$$
note that 
$$
\|P_n^{\alpha,\beta}\|^2_{L^2(\dd\rho_{\alpha,\beta})} =
\frac{2^{\alpha+\beta+1}\Gamma(n+\alpha+1)\Gamma(n+\beta+1)}{(2n+\alpha+\beta+1)\Gamma(n+\alpha+\beta+1)\Gamma(n+1)},
$$
where in case $n=0$ and $\alpha+\beta=-1$ the product $(2n+\alpha+\beta+1)\Gamma(n+\alpha+\beta+1)$ is understood as
$\Gamma(\alpha+\beta+2)$.
The oscillating series defining $K_t^{J,\alpha,\beta}(x,y)$ cannot be computed, nevertheless sharp estimates of the Jacobi
heat kernel are available.

For $\theta,\varphi \in [0,\pi]$ and $t > 0$ define
$$
Z_t^{\alpha,\beta}(\theta,\varphi) := [t+\theta\varphi]^{-\alpha-1/2} [t+ (\pi-\theta)(\pi-\varphi)]^{-\beta-1/2}
	\frac{1}{\sqrt{t}} \exp\bigg( - \frac{(\theta-\varphi)^2}{4t}\bigg).
$$
The following sharp bounds were proved in a series of papers by Sj\"ogren, Szarek and the author
\cite{NSS1,NSS2,NSS3} (see also a related recent article \cite{NPS}).
\begin{thm}[{\cite{NSS1,NSS2,NSS3}}] \label{thm:jac_est}
Let $\alpha,\beta > -1$ and $T > 0$ be fixed. There exists a constant $c=c(\alpha,\beta,T)> 1$ such that
$$
c^{-1}\cdot Z_t^{\alpha,\beta}(\theta,\varphi) \le K_t^{J,\alpha,\beta}(\cos\theta,\cos\varphi)
	\le c\cdot Z_t^{\alpha,\beta}(\theta,\varphi)
$$
for $\theta,\varphi \in [0,\pi]$ and $0 < t \le T$, and\,\footnote{
\,The large time bounds are a direct consequence of the short time estimates, see e.g.\ \cite[p.\,233]{NoSj}.
} 
$$
c^{-1} \le K_t^{J,\alpha,\beta}(x,y) \le c
$$
for $x,y \in [-1,1]$ and $t \ge T$.
\end{thm}

\section{Mixed Hermite--Laguerre--Jacobi setting} \label{sec:mix}

Let $d \ge 1$ and $D = \{1,\ldots,d\}$.
Let $D^H$, $D^L$ and $D^J$ be mutually disjoint subsets of $D$ such that
$$
D = D^H \sqcup D^L \sqcup D^J;
$$
some of these subsets may be empty. Denote
$$
d^{H} := \big|D^{H}\big|, \qquad d^{L} := \big|D^{L}\big|, \qquad d^{J} := \big|D^{J}\big|,
$$
and consider multi-indices
$$
\alpha = (\alpha_i)_{i \in D^L \sqcup D^J} \in (-1,\infty)^{d^L + d^J} \qquad \textrm{and} \qquad
\beta = (\beta_i)_{i \in D^J} \in (-1,\infty)^{d^J}.
$$
In what follows we assume that $D^H$, $D^L$, $D^J$, $\alpha$ and $\beta$ are fixed.
In order to avoid baroque notation, dependence of various objects on these sets and parameters will not be explicitly indicated.

Bring in the space
\begin{equation} \label{space}
\Omega := \omega_1 \times \ldots \times \omega_d, \qquad \textrm{where} \quad \omega_i :=
	\begin{cases}
		\mathbb{R} & \textrm{if}\;\; i \in D^H \\
		\mathbb{R}_+ & \textrm{if}\;\; i \in D^L \\
		[-1,1] & \textrm{if} \;\; i \in D^J
	\end{cases},
	\qquad i = 1,\ldots,d,
\end{equation}
and a product measure on $\Omega$
\begin{equation} \label{meas}
\dd M(x) := \prod_{i=1}^d \dd m_i(x_i),
\end{equation}
where
$$
\dd m_i(x_i) :=
	\begin{cases}
		\dd \gamma(x_i) & \textrm{if}\;\; i \in D^H \\
		\dd \mu_{\alpha_i}(x_i) & \textrm{if}\;\; i \in D^L \\
		\dd \rho_{\alpha_i,\beta_i}(x_i) & \textrm{if} \;\; i \in D^J
	\end{cases},
	\qquad i = 1,\ldots,d.
$$
Note that $M$ is finite.

Next, define mixed Hermite--Laguerre--Jacobi polynomials
$$
\Phi_n := \prod_{i=1}^d \varphi_{n_i}, \qquad n \in \mathbb{N}^d,
$$
where
$$
\varphi_{n_i} :=
	\begin{cases}
		 H_{n_i} & \textrm{if}\;\; i \in D^H \\
		 L_{n_i}^{\alpha_i} & \textrm{if}\;\; i \in D^L \\
		 P_{n_i}^{\alpha_i,\beta_i} & \textrm{if} \;\; i \in D^J
	\end{cases},
	\qquad i = 1,\ldots,d.
$$
Then the system $\{\Phi_n : n \in \mathbb{N}^d\}$ forms an orthogonal basis in $L^2(\Omega,\dd M)$.
Note that
$$
\|\Phi_n\|_{L^2(\dd M)} = \prod_{i \in D^H} \|H_{n_i}\|_{L^2(\dd \gamma)}
\prod_{j \in D^L} \|L_{n_j}^{\alpha_j}\|_{L^2(\dd\mu_{\alpha_j})}
\prod_{k \in D^J} \|P_{n_k}^{\alpha_k,\beta_k}\|_{L^2(\dd\rho_{\alpha_k,\beta_k})}.
$$

The mixed polynomials $\Phi_n$ are eigenfunctions of the ``Laplacian''
$$
\DD := \sum_{i=1}^d \DD_i, \qquad \textrm{where} \quad \DD_i :=
	\begin{cases}
		\DD_H & \textrm{if}\;\; i \in D^H \\
		\DD_{L,\alpha_i} & \textrm{if}\;\; i \in D^L \\
		\DD_{J,\alpha_i,\beta_i} & \textrm{if} \;\; i \in D^J
	\end{cases};
$$
here each one-dimensional differential operator $\DD_i$ is understood to act on the $i$th coordinate.
More precisely, one has
$$
\DD \Phi_n = \Lambda_n \Phi_n, \qquad n \in \mathbb{N}^d,
$$
with the eigenvalues
$$
\Lambda_n := \sum_{i \in D^H} n_i + 2 \sum_{j \in D^L} n_j +  \sum_{k \in D^J} n_k(n_k+\alpha_k+\beta_k+1).
$$

The operator $\DD$ is symmetric and non-negative in $L^2(\Omega,\dd M)$ when considered initially on the subspace of polynomials.
Its natural and unique self-adjoint extension in $L^2(\Omega,\dd M)$
(denoted by the same symbol) is non-negative and given by
$$
\DD f = \sum_{n \in \mathbb{N}^d} \Lambda_n \frac{\langle f,\Phi_n\rangle_{\dd M} \; \Phi_n}{\|\Phi_n\|^2_{L^2(\dd M)}}
$$
on the domain
$$
\domain \DD := \Bigg\{ f \in L^2(\dd M) : \sum_{n \in \mathbb{N}^d} \bigg| \Lambda_n\;
	\frac{\langle f, \Phi_n\rangle_{\dd M}}{\|\Phi_n\|_{L^2(\dd M)}} \bigg|^2 < \infty \Bigg\}.
$$

The corresponding semigroup
$$
\mathcal{T}_t := \exp(-t \DD), \qquad t \ge 0,
$$
which we call the hybrid Hermite--Laguerre--Jacobi semigroup, has an integral representation on $L^2(\dd M)$
$$
\mathcal{T}_tf(x) = \int_{\Omega} \mathcal{K}_t(x,y) f(y)\, \dd M(y), \qquad x \in \Omega, \quad t > 0,
$$
which is valid\,\footnote{
\,More precisely, the integral converges absolutely for all $f \in L^1(\dd M)$, $x \in \Omega$ and $t>0$ and thus
defines $\mathcal{T}_t$ on $L^1(\dd M)$. Notice that $L^p(\dd M) \subset L^1(\dd M)$, $p > 1$, since the measure is finite.
}
for $f \in L^1(\dd M)$. Here the kernel is given by
$$
\mathcal{K}_t(x,y) := \sum_{n \in \mathbb{N}^d} e^{-t\Lambda_n} \, \frac{\Phi_n(x) \Phi_n(y)}{\|\Phi_n\|^2_{L^2(\dd M)}}
	= \prod_{i=1}^d \mathcal{K}_t^i(x_i,y_i), \qquad x,y \in \Omega, \quad t > 0,
$$
where
$$
\mathcal{K}_t^i(x_i,y_i) := 
	\begin{cases}
		K_t^H(x_i,y_i) & \textrm{if}\;\; i \in D^H \\
		K_t^{L,\alpha_i}(x_i,y_i) & \textrm{if}\;\; i \in D^L \\
		K_t^{J,\alpha_i,\beta_i}(x_i,y_i) & \textrm{if} \;\; i \in D^J
	\end{cases},
	\qquad i = 1,\ldots,d.
$$
Note that $\mathcal{K}_t(x,y)$ is strictly positive.
From the product structure and well known properties of the Hermite, Laguerre and Jacobi semigroups it follows that
$\{\mathcal{T}_t\}_{t \ge 0}$ is a symmetric diffusion semigroup in the sense of \cite[Chapter~III]{St}.

\section{Main result} \label{sec:main}

The principal object of our study is the maximal operator of the hybrid Hermite--Laguerre-Jacobi semigroup
$$
\mathcal{M}f := \sup_{t > 0} |\mathcal{T}_t f|.
$$
By virtue of Stein's maximal theorem \cite[p.\,73]{St}, we see that $\mathcal{M}$ is bounded on $L^p(\dd M)$ for $p > 1$.
This could also be concluded directly from $L^p$-boundedness of the Hermite, Laguerre and Jacobi semigroups maximal operators.

On the other hand, none of the semigroup maximal operators just mentioned is bounded on $L^1$, and the natural substitute
result for $\mathcal{M}$ is the weak type $(1,1)$ estimate
$$
M\Big(\big\{x \in \Omega : \mathcal{M}f(x) > \lambda\big\}\Big)
	\lesssim \frac{1}{\lambda}\|f\|_{L^1(\dd M)}, \qquad \lambda > 0, \quad f \in L^1(\dd M),
$$
which is equivalent to boundedness of $\mathcal{M}$ 
from $L^1(\dd M)$ to weak $L^1(\dd M)$, i.e.\ to the Lorentz space $L^{1,\infty}(\dd M)$. Such result for $\mathcal{M}$ does not follow
from any known general theorem. Neither it is a direct consequence of weak type $(1,1)$ bounds for (multi-dimensional)
Hermite, Laguerre and Jacobi semigroups maximal operators.

The main result of this paper reads as follows. We assume that $d$, $D_H$, $D_L$, $D_J$, $\alpha$ and $\beta$
are as described in the beginning of Section \ref{sec:mix}.
\begin{bigthm} \label{thm:main}
The maximal operator $\mathcal{M}$ of the hybrid Hermite--Laguerre--Jacobi semigroup is bounded
from $L^1(\Omega,\dd M)$ to weak $L^1(\Omega,\dd M)$.
\end{bigthm}

Combined with standard arguments\,\footnote{
\,See e.g.\ \cite[Chapter 2, Section 2]{Duo} and observe that the system $\{\Phi_n : n \in \mathbb{N}^d\}$
is linearly dense in $L^1(\Omega,\dd M)$.
},
Theorem~\ref{thm:main} implies almost everywhere boundary convergence for the hybrid Hermite--Laguerre--Jacobi semigroup
applied to $L^1$ functions.
\begin{corollary} \label{cor:main}
For each $f \in L^1(\Omega,\dd M)$,
$$
	\lim_{t \to 0^+} \mathcal{T}_t f = f \quad \emph{a.e.}
$$
\end{corollary}

The remaining part of the paper is devoted to the proof of Theorem \ref{thm:main}.

\section{Proof of Theorem {\ref{thm:main}}} \label{sec:proof}

We start with several lemmas that will be needed in the proof.

\subsection{Technical preparation}

Recall that $D = \{1,\ldots,d\}$ for some $d \ge 1$. Consider a partition
$$
D = D' \sqcup D''
$$
and denote $d' = |D'|$, $d'' = |D''|$, so that $d'+d'' = d$.
Our product space $\Omega = \omega_1 \times \ldots \times \omega_d$ can be factorized as
$$
\Omega = \Omega' \times \Omega'', \qquad \textrm{where} \quad
	\Omega' = \prod_{i \in D'} \omega_i, \quad \Omega'' = \prod_{i \in D''} \omega_i.
$$
Then points in $\Omega$ can be written as $x = (x',x'') \in \Omega' \times \Omega''$.
The measure factorizes in an analogous way,
$$
\dd M(x) = \dd M'(x') \, \dd M''(x''),
$$
where
$$
\dd M'(x') = \prod_{i \in D'} \dd m_i(x_i), \qquad \dd M''(x'') = \prod_{i \in D''} \dd m_i(x_i).
$$

Next, consider an operator $T$ acting on $L^1(\Omega,\dd M)$ and operators $T'$ and $T''$ acting on
$L^1(\Omega',\dd M')$ and $L^1(\Omega'',\dd M'')$, respectively, such that
$$
|Tf| \le (T' \circ T'')f, \qquad f \in L^1(\Omega,\dd M);
$$
in the composition $T' \circ T''$ the operator $T'$ acts only in coordinates $x'$ and $T''$ acts only in coordinates $x''$
(here we tacitly assume that the composition makes sense).

We now formulate three simple lemmas. Most of these results is well known, at least as a folklore.
Some of the lemmas are stated in slightly greater generality than actually needed, just to give wider perspective for
understanding the arguments justifying our main result.
\begin{lemma} \label{lem:1}
Assume that $T''$ is bounded on $L^1(\Omega'',\dd M'')$.
\begin{itemize}
\item[(a)]
	If $T'$ is bounded on $L^1(\Omega',\dd M')$, then $T$ is bounded on $L^1(\Omega,\dd M)$.
\item[(b)]
	If $T'$ is bounded from $L^1(\Omega',\dd M')$ to weak $L^1(\Omega',\dd M')$,
		then $T$ is bounded from $L^1(\Omega,\dd M)$ to weak $L^1(\Omega,\dd M)$.
\end{itemize}
\end{lemma}

\begin{proof}
Item (a) is obvious. To verify (b), apply first the weak type $(1,1)$ estimate in the $x'$ variables and then strong type $(1,1)$
bound in the $x''$ variables.
\end{proof}

Consider a bimeasurable bijection $\Psi\colon \Omega \mapsto \widetilde{\Omega}$ and let $\Psi \# M$
be the push-forward measure on $\widetilde{\Omega}$ induced by $\Psi$.
The next result is a simple observation, which we nevertheless formulate as a lemma for easy further reference.
\begin{lemma} \label{lem:2}
$T$ is bounded from $L^1(\Omega,\dd M)$ to weak $L^1(\Omega,\dd M)$ if and only if the operator $T_{\Psi}$ defined by
$$
(T_{\Psi}f) \circ \Psi = T(f \circ \Psi)
$$
is bounded from $L^1(\widetilde{\Omega},\dd (\Psi \# M))$ to weak $L^1(\widetilde{\Omega},\dd (\Psi \# M))$.
\end{lemma}
Note that when $T$ is, for instance, of the form $Tf(x) = \sup_t \int_{\Omega} K_t(x,y) |f(y)|\, \dd M(y)$
and $\dd M(y) = \mathfrak{m}(y)\, \dd y$
with differentiable density $\mathfrak{m}$, then $T_{\Psi}$ from Lemma~\ref{lem:2} is of the form
$T_{\Psi}F(x) = \sup_t \int_{\widetilde{\Omega}} K_t(\Psi^{-1}(x),\Psi^{-1}(y)) |F(y)|\, \dd (\Psi\# M)(y)$, where
$\dd (\Psi \# M)(y) = \mathfrak{m} \circ \Psi^{-1}(y) |(\Psi^{-1})'(y)|\, \dd y$.

Let now $\Omega'$ and $M'$ be invariant under reflections with respect to (vector) hyperplanes perpendicular to coordinate
axes in $\mathbb{R}^{d'}$. Denote
$$
\Omega'_R = \Omega' \cap \mathbb{R}^{d'}_{+}, \qquad M'_R = M'\big|_{\Omega'_R}
$$
and
$$
\Omega_R = \Omega'_R \times \Omega'', \qquad M_R = M'_R \times M''.
$$
Let $T_R$ be an operator on $L^1(\Omega_R,\dd M_R)$ given by
$$
T_R f(x) = Tf^{e}(x), \qquad x \in \Omega_R,
$$
where $f^e$ is the unique extension of $f$ to $\Omega$ that is even with respect to all the $x'$ coordinates.
For $\eta = (\eta_i)_{i \in D'} \in \{-1,1\}^{d'}$ let $\Psi_{\eta}\colon \Omega \mapsto \Omega$ be defined by 
$$
\big( \Psi_{\eta}(x)\big)_i = 
	\begin{cases}
		\eta_i x_i & \textrm{if}\;\; i \in D' \\
		x_i & \textrm{if}\;\; i \in D''
	\end{cases},
	\qquad i=1,\ldots,d, \quad x \in \Omega.
$$

\begin{lemma} \label{lem:3}
Assume that $T$, considered on $L^1(\Omega,\dd M)$, satisfies the following conditions.
\begin{itemize}
\item[(i)] $T$ is sublinear.
\item[(ii)] $|Tf| \le T|f|$.
\item[(iii)] $Tf_1 \le Tf_2$ whenever $0 \le f_1 \le f_2$.
\item[(iv)] $T$ is $(\Omega',\dd M')$-symmetric in the sense that 
$$
(Tf) \circ \Psi_{\xi} = Tf \qquad \textrm{for all} \;\; \xi \in \{-1,1\}^{d'}.
$$
\end{itemize}
Then $T$ is bounded from $L^1(\Omega,\dd M)$ to weak $L^1(\Omega,\dd M)$ if and only if $T_R$ is bounded
from $L^1(\Omega_R,\dd M_R)$ to weak $L^1(\Omega_R,\dd M_R)$.
\end{lemma}

To prove the lemma it is convenient to introduce some additional notation and terminology.
Let $\epsilon = (\epsilon_i)_{i \in D'} \in \{0,1\}^{d'}$. We say that a function $f\colon \Omega \mapsto \mathbb{C}$ is
$\epsilon$-symmetric if for each $i \in D'$, $f$ is either even or odd with respect to the $i$th coordinate according
to whether $\epsilon_i = 0$ or $\epsilon_i = 1$, respectively. Define $\epsilon$-symmetric component of an arbitrary
$f\colon \Omega \mapsto \mathbb{C}$ by
$$
f_{\epsilon} = \frac{1}{2^{d'}} \sum_{\xi \in \{-1,1\}^{d'}} \xi^{\epsilon} f \circ \Psi_{\xi}, \qquad \textrm{where} \;\;
	\xi^{\epsilon} = \prod_{i \in D'} \xi_i^{\epsilon_i}.
$$
Then
$$
f = \sum_{\epsilon \in \{0,1\}^{d'}} f_{\epsilon}.
$$
For an $\epsilon$-symmetric $f$, denote by $f_R$ its restriction to $\Omega_R$.
Recall that for a function $f$ acting on $\Omega_R$, 
by $f^e$ we mean its $(0,\ldots,0)$-symmetric extension to $\Omega$.

\begin{proof}[{Proof of Lemma \ref{lem:3}}]
First, assume that $T_R$ is satisfies the weak type $(1,1)$ bound.
Then, in view of our assumptions (i)--(iv) on $T$, for $f \in L^1(\Omega,\dd M)$ and $\lambda > 0$ we have
\begin{align*}
M\Big( \big\{ x \in \Omega : |Tf(x)| > \lambda \big\} \Big) & = 
	M\bigg( \bigg\{ x \in \Omega : \bigg|T\Bigg( \sum_{\epsilon \in \{0,1\}^{d'}} f_{\epsilon} \Bigg)(x)\bigg| > \lambda \bigg\}\bigg) \\
& \overset{(ii)}{\le} 
M\bigg( \bigg\{ x \in \Omega :  T \Bigg|\sum_{\epsilon \in \{0,1\}^{d'}} f_{\epsilon}\Bigg| (x)> \lambda \bigg\}\bigg) \\
& \overset{(iii)}{\le} 
M\bigg( \bigg\{ x \in \Omega :  T \Bigg(\sum_{\epsilon \in \{0,1\}^{d'}} |f_{\epsilon}|\Bigg) (x)> \lambda \bigg\}\bigg) \\
& \overset{(i)}{\le}
M\bigg( \bigg\{ x \in \Omega :  \sum_{\epsilon \in \{0,1\}^{d'}} T|f_{\epsilon}|(x)> \lambda \bigg\}\bigg) \\
& \le  \sum_{\epsilon \in \{0,1\}^{d'}} M\Big( \big\{ x \in \Omega : T|f_{\epsilon}|(x)> \lambda/2^{d'} \big\}\Big) \\
& = \sum_{\epsilon \in \{0,1\}^{d'}} \sum_{\xi \in \{-1,1\}^{d'}} 
		M_R\Big( \big\{ x \in \Omega_R : (T|f_{\epsilon}|)\circ \Psi_{\xi}(x) > \lambda/2^{d'}\big\}\Big) \\
& \overset{(iv)}{=}
\sum_{\epsilon \in \{0,1\}^{d'}} 2^{d'}
		M_R\Big( \big\{ x \in \Omega_R : T|f_{\epsilon}|(x) > \lambda/2^{d'}\big\}\Big) \\
& = \sum_{\epsilon \in \{0,1\}^{d'}} 2^{d'}
		M_R\Big( \Big\{ x \in \Omega_R :  T_R\big(|f_{\epsilon}|_R \big)(x) > \lambda/2^{d'}\Big\}\Big) \\
& \lesssim \sum_{\epsilon \in \{0,1\}^{d'}}
		\frac{4^{d'}}{\lambda} \big\||f_{\epsilon}|_R \big\|_{L^1(\Omega_R,\dd M_R)} 
= \sum_{\epsilon \in \{0,1\}^{d'}}
		\frac{2^{d'}}{\lambda} \big\|f_{\epsilon}\big\|_{L^1(\Omega,\dd M)} \\
& \le \frac{4^{d'}}{\lambda} \|f\|_{L^1(\Omega,\dd M)},
\end{align*}
where we also used the facts that $|f_{\epsilon}|$ is even with respect to all the $x'$ coordinates and that
$\|f_{\epsilon}\|_{L^1(\Omega,\dd M)} \le \|f\|_{L^1(\Omega,\dd M)}$. Thus $T$ is weak type $(1,1)$.

To verify the other implication,
assume now that $T$ satisfies the weak type $(1,1)$ estimate.
Then, for $f \in L^1(\Omega_R,\dd M_R)$ and $\lambda > 0$, we can write
\begin{align*}
M_R\Big( \big\{ x \in \Omega_R : |T_R f(x)| > \lambda \big\}\Big)
& = M_R\Big( \big\{ x \in \Omega_R : \big|T f^{e}(x)\big| > \lambda \big\}\Big) \\
& \le M\Big( \big\{ x \in \Omega : \big|Tf^{e}(x)\big| > \lambda \big\} \Big) \\
& \lesssim \frac{1}{\lambda} \big\|f^{e}\big\|_{L^1(\Omega,\dd M)} \\
& = \frac{2^{d'}}{\lambda} \|f\|_{L^1(\Omega_R,\dd M_R)}.
\end{align*}

The conclusion follows.
\end{proof}

\subsection{Reduction to the case $D_H = \emptyset$}
This step is based on relating, in a one-dimensional situation, the integral kernel $K_t^H(x,y)$
of the Ornstein--Uhlenbeck semigroup with the kernel $K_t^{L,\alpha}(x,y)$ of the Laguerre semigroup.
Explicit formulas for these kernels are given in Sections \ref{ssec:H} and \ref{ssec:L}, respectively. 

Observe that $K_t^H(x,y)$ can be decomposed as
$$
K_t^H(x,y) = K_t^{H,e}(x,y) + K_t^{H,o}(x,y), \qquad x,y \in \mathbb{R},
$$
where
\begin{align*}
K_t^{H,e}(x,y) & := \sqrt{\frac{e^t}{2\pi \sinh t}} \exp\bigg( - \frac{e^{-t}(x^2+y^2)}{2\sinh t}\bigg)
	\cosh \bigg(\frac{|xy|}{\sinh t}\bigg), \\
K_t^{H,o}(x,y) & := \sqrt{\frac{e^t}{2\pi \sinh t}} \exp\bigg( - \frac{e^{-t}(x^2+y^2)}{2\sinh t}\bigg)
	\sgn(xy)\, \sinh \bigg(\frac{|xy|}{\sinh t}\bigg).
\end{align*}
Clearly,
\begin{equation} \label{HLmaj}
K_t^H(x,y) \le K_t^{H,e}(x,y) + \big|K_t^{H,o}(x,y)\big| \le 2 K_t^{H,e}(x,y), \qquad x,y \in \mathbb{R}.
\end{equation}

On the other hand, using the identity (cf.\ \cite[Chapter~5, (5.8.5)]{Leb})
$$
I_{-1/2}(z) = \sqrt{\frac{2}{\pi z}} \cosh z
$$
we see that
\begin{equation*}
K_t^{H,e}(x,y) = \frac{e^{t/2}}{2\sinh t} \exp\bigg( - \frac{e^{-t}(x^2+y^2)}{2\sinh t}\bigg) \sqrt{|xy|}\,
	I_{-1/2}\bigg( \frac{|xy|}{\sinh t}\bigg),
\end{equation*}
with the limiting interpretation in case $xy = 0$. Thus
\begin{equation} \label{HLcon}
K_t^{H,e}(x,y) = K^{L,-1/2}_t\big(x^2,y^2\big).
\end{equation}

Now, recall that the integral kernel of the hybrid semigroup is given by
\begin{equation}  \label{ker}
\mathcal{K}_t(x,y) = \prod_{i \in D^H} K_t^H(x_i,y_i) \prod_{j \in D^L} K_t^{L,\alpha_j}(x_j,y_j)
	\prod_{k \in D^J} K_t^{J,\alpha_k,\beta_k}(x_k,y_k),
\end{equation}
and the corresponding measure space is $(\Omega,\dd M) = (\omega_1\times\ldots\times\omega_d,\dd m_1\ldots\dd m_d)$,
see \eqref{space} and \eqref{meas}.

Assume that $D^H \neq \emptyset$.
Because of \eqref{HLmaj}, it is sufficient to prove Theorem \ref{thm:main} with $K_t^H(x_i,y_i)$ in \eqref{ker} replaced by
$K_t^{H,e}(x_i,y_i)$; the multiplicative constant $2^{d^H}$ is of course irrelevant for our purpose and thus omitted
(we will skip analogous comments in the sequel).

Then, in view of the reflection symmetries of $K_t^{H,e}(x_i,y_i)$ and Lemma \ref{lem:3}
(specified to $\Omega' = \prod_{i \in D^H} \omega_i$),
one can further replace factors $\omega_i = \mathbb{R}$ corresponding to $i \in D^H$ in the product defining the space by
$\omega_i = \mathbb{R}_+$ and restrict the measure accordingly.

Finally, by \eqref{HLcon} and Lemma \ref{lem:2} (specified to $\Psi$ given by $(\Psi(x))_i = x_i^2$ if $i \in D^H$,
otherwise $(\Psi(x))_i = x_i$), one can further replace the factors $K_t^{H,e}(x_i,y_i)$ by $K_t^{L,-1/2}(x_i,y_i)$, setting
simultaneously $\dd m_i(x_i) = \dd \mu_{-1/2}(x_i)$ for $i \in D^H$.

Thus we reduced proving Theorem \ref{thm:main} to the case $D^H = \emptyset$.

\subsection{Further reduction to the case $D_J = D_H = \emptyset$}
From now we assume that $D_H = \emptyset$ and hence the kernel of the hybrid semigroup is
\begin{equation} \label{zz}
\mathcal{K}_t(x,y) = \prod_{j \in D^L} K_t^{L,\alpha_j}(x_j,y_j) \prod_{k \in D^J} K_t^{J,\alpha_k,\beta_k}(x_k,y_k).
\end{equation}
We shall establish a control of the corresponding maximal operator $\mathcal{M}$ by a sum of finite number of auxiliary
maximal operators. Then we will show that each component of the sum satisfies the weak type $(1,1)$ estimate
with respect to the measure space $(\Omega,\dd M)$, provided that the multi-dimensional Laguerre semigroup maximal operator
is of weak type $(1,1)$. This is obvious when initially $D_J=\emptyset$, so in what follows we assume that $D_J \neq \emptyset$.

For $k \in D^J$ and $x_k,y_k \in [-1,1]$ denote
$$
\theta_k = \theta_k(x_k) := \arccos x_k, \qquad \varphi_k = \varphi_k(y_k) := \arccos y_k.
$$
Observe that for points $x,y \in \Omega$ and $t > 0$
\begin{align}
1 & \le \chi_{\{t > 1\}} + \chi_{\{t \le 1\}} \prod_{k \in D^J} \Big( \chi_{\{\theta_k,\varphi_k \le 3\pi/4\}}
	+ \chi_{\{\theta_k,\varphi_k \ge \pi/4\}} + \chi_{\{|\theta_k-\varphi_k| \ge \pi/2\}} \Big) \nonumber \\
& = \chi_{\{t > 1\}} + \sum_{A \sqcup B \sqcup C = D^J } \chi_{\{t \le 1\}} \prod_{i \in A} \chi_{\{\theta_i,\varphi_i \le 3\pi/4\}}
	\prod_{j \in B} \chi_{\{\theta_j,\varphi_j \ge \pi/4\}} \prod_{k \in C}\chi_{\{|\theta_k-\varphi_k| \ge \pi/2\}}, \label{decomp}
\end{align}
where the last summation runs over all weak\,\footnote{
\,Here \emph{weak} refers to partitions allowing empty sets.
}
ordered partitions of $D^J$ into subsets $A,B,C$.
The expression in \eqref{decomp} is a sum of $1+3^{d^J}$ terms. Multiplying our kernel $\mathcal{K}_t(x,y)$ in \eqref{zz}
by each of these terms separately leads to $1+3^{d^J}$ auxiliary kernels, call them $\mathcal{K}_t^{\infty}(x,y)$
and $\mathcal{K}_t^{A,B,C}(x,y)$, respectively.
Then
$$
\mathcal{K}_t(x,y) \le \mathcal{K}_t^{\infty}(x,y) + \sum_{A \sqcup B \sqcup C = D^J} \mathcal{K}_t^{A,B,C}(x,y).
$$
Consequently,
$$
\mathcal{M}f \le \mathcal{M}_{\infty}f + \sum_{A \sqcup B \sqcup C = D^J} \mathcal{M}_{A,B,C}f,
$$
where the auxiliary maximal operators are defined by
\begin{align*}
\mathcal{M}_{\infty}f(x) & = \sup_{t > 1} \int_{\Omega} \mathcal{K}_t^{\infty}(x,y) |f(y)|\, \dd M(y), \\
\mathcal{M}_{A,B,C}f(x) & = \sup_{0 < t \le 1} \int_{\Omega} \mathcal{K}_t^{A,B,C}(x,y) |f(y)|\, \dd M(y).
\end{align*}

We now analyze $\mathcal{M}_{\infty}$ and $\mathcal{M}_{A,B,C}$. It will be convenient to use the notation
$$
\Omega' = \prod_{j \in D^L} \omega_j, \qquad \dd M' = \prod_{j \in D^L} \dd m_j, \qquad
\Omega'' = \prod_{k \in D^J} \omega_k, \qquad \dd M'' = \prod_{k \in D^J} \dd m_k
$$
and
$$
\Omega_{\Diamond}'' = \prod_{k \in \Diamond} \omega_k, \qquad \dd M_{\Diamond}'' = \prod_{k \in \Diamond} \dd m_k, \qquad
	\Diamond = A,B,C.
$$

\subsubsection*{Treatment of $\mathcal{M}_{\infty}$}
We have the following.
\begin{lemma} \label{lem:inf}
The maximal operator
$$
f \mapsto \sup_{t \ge 1} \int_{\Omega''} \prod_{k \in D^J} K_t^{J,\alpha_k,\beta_k}(x_k,y_k) |f(y)| \, \dd M''(y)
$$
is bounded on $L^1(\Omega'',\dd M'')$.
\end{lemma}

\begin{proof}
This is an immediate consequence of Theorem \ref{thm:jac_est}, since for $t \ge 1$ the one-dimensional kernels
$K_t^{J,\alpha_k,\beta_k}(x_k,y_k)$ are dominated by constants independent of $t$ and $x_k,y_k$.
\end{proof}

The operator $\mathcal{M}_{\infty}$ is controlled by a composition of the Laguerre maximal operator
$$
f \mapsto \sup_{t>0} \int_{\Omega'} \bigg[ \prod_{j \in D^L} K_t^{L,\alpha_j}(x_j,y_j)\bigg] |f(y)|\, \dd M'(y)
$$
acting on the $D^L$ coordinates and the maximal operator from Lemma~\ref{lem:inf} acting on the $D^J$ coordinates.
Thus, in view of the boundedness claimed in Lemma~\ref{lem:inf},
we are in the situation of Lemma~\ref{lem:1}(b). It follows that $\mathcal{M}_{\infty}$ is of weak type $(1,1)$
if only the Laguerre maximal operator possesses analogous property.

\subsubsection*{Treatment of $\mathcal{M}_{A,B,C}$, reduction to the case $C = \emptyset$}
We assume $C \neq \emptyset$, otherwise there is nothing to do in this step.
\begin{lemma} \label{lem:far}
The maximal operator
$$
f \mapsto \sup_{0 < t \le 1} \int_{\Omega_{C}''}
	\bigg[ \prod_{k \in C} \chi_{\{|\theta_k-\varphi_k|\ge \pi/2\}} K_t^{J,\alpha_k,\beta_k}(x_k,y_k)\bigg]
	|f(y)| \, \dd M_{C}''(y)
$$
is bounded on $L^1(\Omega_{C}'',\dd M_{C}'')$.
\end{lemma}

\begin{proof}
The argument is essentially the same as in the proof of Lemma~\ref{lem:inf}.
Let $k \in C$. Assume that $0 < t \le 1$ and $|\theta_k-\varphi_k| \ge \pi/2$.
Then, by Theorem~\ref{thm:jac_est},
\begin{align*}
K_t^{J,\alpha_k,\beta_k}(x_k,y_k) & \lesssim [t+\theta_k\varphi_k]^{-\alpha_k-1/2} [t+(\pi-\theta_k)(\pi-\varphi_k)]^{-\beta_k-1/2}
	\frac{1}{\sqrt{t}}\exp\bigg(-\frac{(\theta_k-\varphi_k)^2}{4t}\bigg) \\
& \lesssim [t+\theta_k\varphi_k]^{-\alpha_k-1} [t+(\pi-\theta_k)(\pi-\varphi_k)]^{-\beta_k-1}
	\frac{1}{\sqrt{t}}\exp\bigg(-\frac{(\pi/2)^2}{4t}\bigg) \\
& \lesssim t^{-\alpha_k-\beta_k-5/2} e^{-\pi^2/(16t)} \lesssim 1.
\end{align*}
This shows that the kernel of our maximal operator is controlled by a constant independent of $t$ and the arguments.
The conclusion follows.
\end{proof}

The operator $\mathcal{M}_{A,B,C}$ is dominated by a composition of the maximal operator
\begin{equation} \label{maxAB}
\begin{split} 
 f \mapsto \sup_{0 < t \le 1} \int_{\Omega' \times \Omega_{A}'' \times \Omega_{B}''} |f(y)| 
	\prod_{j \in D^L} K_t^{L,\alpha_j}(x_j,y_j)\, \dd m_j(y_j)  \\ 
\times \prod_{k \in A} \chi_{\{\theta_k,\varphi_k \le 3\pi/4\}}
		K_t^{J,\alpha_k,\beta_k}(x_k,y_k)\, \dd m_k(y_k) \\
\times \prod_{l \in B} \chi_{\{\theta_l,\varphi_l \ge \pi/4\}}
		K_t^{J,\alpha_l,\beta_l}(x_l,y_l)\, \dd m_l(y_l) 
\end{split}
\end{equation}
acting on the $D^L \sqcup A \sqcup B$ coordinates and the maximal operator from Lemma~\ref{lem:far}.
Therefore, taking into account the boundedness from Lemma~\ref{lem:far}, it follows by Lemma~\ref{lem:1}(b) that $\mathcal{M}_{A,B,C}$
is of weak type $(1,1)$ if only $\mathcal{M}_{A,B,\emptyset}$ satisfies weak type $(1,1)$ bound.

\subsubsection*{Treatment of $\mathcal{M}_{A,B,C}$, further reduction to the case $B = C = \emptyset$}
In this step we assume that $C = \emptyset$ and $B \neq \emptyset$.
Consider the maximal operator defined in \eqref{maxAB} (i.e.\ $\mathcal{M}_{A,B,\emptyset}$) and the bijection
$\Psi\colon \Omega' \times \Omega_{A}'' \times \Omega_{B}'' \mapsto  \Omega' \times \Omega_{A}'' \times \Omega_{B}''$
given by $(\Psi(x))_i = -x_i$ if $i \in B$, otherwise $(\Psi(x))_i = x_i$.
Notice that $\Psi^{-1}=\Psi$.
Under this bijection applied in both $x$ and $y$ variables
the first two products in \eqref{maxAB} remain invariant, whereas the last one is transformed to
\begin{equation} \label{maxZ}
\prod_{l \in B} \chi_{\{\theta_l,\varphi_l \le 3\pi/4\}} K_t^{J,\beta_l,\alpha_l}(x_l,y_l)\, \dd \rho_{\beta_l,\alpha_l}(y_l);
\end{equation}
here we used the fact that the Jacobi polynomials satisfy (cf.\ \cite[Chapter IV, (4.1.3)]{Sz})
$$
P_n^{\alpha,\beta}(-z) = (-1)^n P_n^{\beta,\alpha}(z), \qquad n \ge 0.
$$

Observe that the expression in \eqref{maxZ} has the same form as the second product in \eqref{maxAB}. Hence
both products can be merged into one product of the same type.
Thus, in view of Lemma \ref{lem:2}, weak type $(1,1)$ bound for the maximal operator in \eqref{maxAB}
holds provided that analogous estimates hold for maximal operators (with any choice of $A$ and the Jacobi parameters) of the form
\begin{equation} \label{maxA}
\begin{split} 
 f \mapsto \sup_{0 < t \le 1} \int_{\Omega' \times \Omega_{A}''} |f(y)| 
	\prod_{j \in D^L} K_t^{L,\alpha_j}(x_j,y_j)\, \dd m_j(y_j)  \qquad \\ 
\times \prod_{k \in A} \chi_{\{\theta_k,\varphi_k \le 3\pi/4\}}
		K_t^{J,\alpha_k,\beta_k}(x_k,y_k)\, \dd m_k(y_k),
\end{split}
\end{equation}
that is for $\mathcal{M}_{A,\emptyset,\emptyset}$.

\subsubsection*{Treatment of $\mathcal{M}_{A,B,C}$, final reduction to the Laguerre case}
In this step we assume that $B=C=\emptyset$ and $A \neq \emptyset$, i.e.\ $\Omega'' = \Omega''_{A}$.

\noindent\textbf{Step 1.}
Let $\widetilde{\Omega}_{A}'' = \prod_{k \in A} [0,\pi]$.
Consider the maximal operator defined in \eqref{maxA} (i.e.\ $\mathcal{M}_{A,\emptyset,\emptyset}$) and the
bijection $\Psi \colon \Omega'\times \Omega_{A}'' \mapsto \Omega' \times \widetilde{\Omega}_{A}''$ given by
$(\Psi(x))_i = \arccos x_i$ if $i \in A$, otherwise $(\Psi(x))_i = x_i$. Under this bijection applied both in $x$ and $y$ variables
the first product in \eqref{maxA} remains invariant, while the second one is transformed
(see the comment succeeding Lemma~\ref{lem:2}) to
$$
\prod_{k \in A} \chi_{\{x_k,y_k \le 3\pi/4\}} K_t^{J,\alpha_k,\beta_k}(\cos x_k, \cos y_k)\, \dd \widetilde{m}_k(y_k),
$$ 
where $\dd \widetilde{m}_k(y_k) := 2^{\alpha_k + \beta_k + 1} [\sin(y_k/2)]^{2\alpha_k+1}[\cos(y_k/2)]^{2\beta_k+1}\, \dd y_k$.
By Lemma \ref{lem:2} we see that weak type $(1,1)$ bound for the maximal operator in \eqref{maxA} is equivalent to weak type
$(1,1)$ estimate for
\begin{equation} \label{maxAtilde}
\begin{split} 
 f \mapsto \sup_{0 < t \le 1} \int_{\Omega' \times \widetilde{\Omega}_{A}''} |f(y)| 
	\prod_{j \in D^L} K_t^{L,\alpha_j}(x_j,y_j)\, \dd m_j(y_j)  \qquad \qquad \quad \\ 
\times \prod_{k \in A} \chi_{\{x_k,y_k \le 3\pi/4\}}
		K_t^{J,\alpha_k,\beta_k}(\cos x_k,\cos y_k)\, \dd \widetilde{m}_k(y_k)
\end{split}
\end{equation}
with respect to the measure space $(\Omega' \times \widetilde{\Omega}_{A}'', \dd M' \prod_{k \in A} \dd \widetilde{m}_k)$.

\noindent\textbf{Step 2.} 
We observe that in \eqref{maxAtilde}, thanks to the characteristic function $\chi_{\{x_k,y_k \le 3\pi/4\}}$, replacing the
measures $\dd \widetilde{m}_k(y_k)$ by comparable measures
\begin{equation} \label{mtt_form}
\dd {\breve{m}}_k(y_k) := y_k^{2\alpha_k+1}\, \dd y_k
\end{equation}
leads to a maximal operator whose weak type $(1,1)$ boundedness with respect to the measure space
$(\Omega' \times \widetilde{\Omega}_{A}'', \dd M' \prod_{k \in A} \dd {\breve{m}}_k)$
is equivalent to weak type $(1,1)$ bound for for the maximal operator in $\eqref{maxAtilde}$.
Moreover, this equivalence is preserved after replacing $\widetilde{\Omega}_{A}''$ with
$$
{{\breve{\Omega}}}''_{A} := \prod_{k \in A} \mathbb{R}_+
$$
and defining $\dd {\breve{m}}_k$ on $\mathbb{R}_+$ by the same formula \eqref{mtt_form}.
That is, we are now concerned with the maximal operator
\begin{equation} \label{maxAtt}
\begin{split} 
 f \mapsto \sup_{0 < t \le 1} \int_{\Omega' \times {\breve{\Omega}}_{A}''} |f(y)| 
	\prod_{j \in D^L} K_t^{L,\alpha_j}(x_j,y_j)\, \dd m_j(y_j)  \qquad \qquad \quad \\ 
\times \prod_{k \in A} \chi_{\{x_k,y_k \le 3\pi/4\}}
		K_t^{J,\alpha_k,\beta_k}(\cos x_k,\cos y_k)\, \dd {\breve{m}}_k(y_k)
\end{split}
\end{equation}
relative to the measure space 
$(\Omega' \times {\breve{\Omega}}_{A}'', \dd M' \prod_{k \in A} \dd {\breve{m}}_k)$.

\noindent\textbf{Step 3.} 
Consider now the bijection
$\Psi\colon \Omega' \times {\breve{\Omega}}_{A}'' \mapsto \Omega' \times {\breve{\Omega}}_{A}'' $
given by $(\Psi(x))_i = x^2_i/2$ for $i \in A$, otherwise $(\Psi(x))_i = x_i$.
Then, with the aid of Lemma~\ref{lem:2}, we see that weak type $(1,1)$ bound for the maximal operator \eqref{maxAtt}
is equivalent to weak type $(1,1)$ estimate for
\begin{equation} \label{maxAttt}
\begin{split} 
 f \mapsto \sup_{0 < t \le 1} \int_{\Omega' \times {\breve{\Omega}}_{A}''} |f(y)| 
	\prod_{j \in D^L} K_t^{L,\alpha_j}(x_j,y_j)\, \dd m_j(y_j) \qquad \qquad \qquad \qquad \;\; \\ 
\times \prod_{k \in A} \chi_{\{x_k,y_k \le 9\pi^2/32\}}
		K_t^{J,\alpha_k,\beta_k}(\cos \sqrt{2 x_k},\cos \sqrt{2 y_k})\, \dd \mu_{\alpha_k}(y_k)
\end{split}
\end{equation}
with respect to the measure space 
$(\Omega' \times {\breve{\Omega}}_{A}'', \dd M' \prod_{k \in A} \dd \mu_{\alpha_k})$.

\noindent\textbf{Step 4.} 
For the final stroke, we need the following.
\begin{lemma} \label{lem:lag_est}
Let $\alpha > -1$. The one-dimensional Laguerre kernel $K_t^{L,\alpha}(x,y)$ satisfies
$$
K_t^{L,\alpha}(x,y) \simeq \big(t+\sqrt{xy}\big)^{-\alpha-1/2} \frac{1}{\sqrt{t}}
	\exp\bigg( - \frac{(\sqrt{x}-\sqrt{y})^2}{2t} \bigg)
$$
uniformly in $0 \le x,y \le 3$ and $0 < t \le 1$.
\end{lemma}

\begin{proof}
Throughout the proof we may assume $x,y \in (0,3]$. Let $t \in (0,1]$.
By the explicit formula \eqref{lag_form} and the Bessel function bounds \eqref{bes_est} we have
\begin{align*}
K_t^{L,\alpha}(x,y) & \simeq \frac{e^{t(\alpha+1)}}{2\sinh t} (xy)^{-\alpha/2} \bigg(\frac{\sqrt{xy}}{\sinh t}\bigg)^{\alpha}
	\bigg( \frac{\sqrt{xy}}{\sinh t} + 1\bigg)^{-\alpha-1/2} \exp\bigg( - \frac{e^{-t}(x+y)-2\sqrt{xy}}{2\sinh t}\bigg) \\
& \simeq \big(t+\sqrt{xy}\big)^{-\alpha-1/2} \frac{1}{\sqrt{t}} \exp\bigg( - \frac{(\sqrt{x}-\sqrt{y})^2}{2t} \bigg) \\
& \qquad \qquad \times	\exp\bigg( \frac{(\sqrt{x}-\sqrt{y})^2}{2t} - \frac{e^{-t}(x+y)-2\sqrt{xy}}{2\sinh t}\bigg),
\end{align*}
the last relation since $\sinh t \simeq t$ under our assumption on $t$.
It is now enough to show that the argument of the last exponential stays bounded for the considered $x,y$ and $t$.

Observe that
\begin{align*}
& \frac{(\sqrt{x}-\sqrt{y})^2}{2t} - \frac{e^{-t}(x+y)-2\sqrt{xy}}{2\sinh t} \\
& \qquad \qquad = \frac{1-e^{-t}}{2t} (x+y) + \frac{1}2 \bigg(\frac{1}t-\frac{1}{\sinh t}\bigg)
	\Big[ \big(\sqrt{x}-\sqrt{y}\big)^2 - (1-e^{-t})(x+y)\Big].
\end{align*}
Since the functions
$$
t \mapsto \frac{1-e^{-t}}{2t} \qquad \textrm{and} \qquad t \mapsto \frac{1}t - \frac{1}{\sinh t}
$$
are bounded on $(0,1]$, we arrive at the desired conclusion.
\end{proof}

By Theorem~\ref{thm:jac_est} and Lemma \ref{lem:lag_est} the second product in \eqref{maxAttt} is controlled by the product
$$
\prod_{k \in A} K_t^{L,\alpha_k}(x_k,y_k)\, \dd \mu_{\alpha_k}(y_k)
$$
of Laguerre kernels and measures, uniformly in $x_k,y_k \ge 0$ and $0 < t \le 1$.
Therefore, the maximal operator in \eqref{maxAttt} is controlled by the multi-dimensional Laguerre maximal operator.
Consequently, weak type $(1,1)$ bound for the latter implies weak type $(1,1)$ estimate for the operator in \eqref{maxAttt}.

Summing up, we have reduced the proof of Theorem~\ref{thm:main} to the case $D^H=D^J=\emptyset$.

\subsection{Epilogue}
For the proof of weak type $(1,1)$ estimate for the multi-dimensional Laguerre semigroup maximal operator we refer
to \cite[Section~2]{NSS}. The reasoning is rather long and somewhat tricky, so we shall not even sketch the details here.
Notice that slightly different time scaling that we use in this paper does not change the Laguerre maximal operator.


\end{document}